\documentclass[review]{elsarticle}
\usepackage{lineno,hyperref}
\usepackage{amsmath}
\usepackage{xcolor}
\usepackage{amssymb}
\usepackage{geometry}
\newtheorem{theorem}{Theorem}
\newtheorem{definition}{Definition}
\newtheorem{lemma}[subsection]{Lemma}
\newdefinition{remark}{Remark}
\newproof{proof}{Proof}
\newproof{pot}{Proof of Theorem \ref{thm2}}

\modulolinenumbers[5]

\journal{Journal of \LaTeX\ Templates}

\begin{document}

\begin{frontmatter}

\title{Inverse problem for recovery of temporal component of source term for multi-term time fractional parabolic equation with nonlocal boundary datum }
\author[label1]{Muhammad Ali}\ead{muhammad.ali.pk.84@gmail.com, m.ali@nu.edu.pk}
\author[label1]{Sara Aziz}\ead{saraazizpk@gmail.com}
\address[label1]{Department of Sciences and Humanities, National University of Computer and Emerging Sciences,\\ Islamabad, Pakistan}

\begin{abstract}
Inverse problem for multi-term fractional parabolic equation in two dimensional space, involving $m+1$ Caputo fractional derivatives in time, is investigated.  Presence of nonlocal boundary conditions leads to a non-self-adjoint spectral problem. A bi-orthogonal system of functions is used to construct the solution that involves double infinite series. Properties of multinomial Mittag-Leffler function and eigenfunctions are used to prove the classical nature of the solution under certain regularity conditions on the given datum.
\end{abstract}

\begin{keyword}
\texttt{Multi-term time fractional differential equation, Smarskii-Ionkin boundary datum, Bi-orthogonal system of functions, Fourier's method, Multinomial Mittag-Leffler function}
\end{keyword}

\end{frontmatter}

\linenumbers
\section{Introduction} \label{intro}

This research article is devoted to the study of inverse problem of extracting the temporal component of the source term $a(t)$ for the following multi-term time fractional  parabolic differential system
\begin{equation}\label{problem}
\bigl({}^cD_{0_+,t}^{\alpha}+\sum_{i=1}^{m}\psi_i
{}^cD_{0_+,t}^{\alpha_i}\bigr)u(x,y,t)+\Delta^2u(x,y,t)=a(t)f(x,y,t),\qquad (x,y,t)\in \Pi,
\end{equation}
alongside the homogeneous boundary conditions
\begin{eqnarray}\label{dbcs}
&u_{x}(0,y,t)=0=u_{xxx}(0,y,t),\;
u_{y}(x,0,t)=0=u_{y}(x,1,t),&\;\nonumber\\
&u_{yyy}(x,0,t)=0=u_{yyy}(x,1,t),&
\end{eqnarray}
 nonlocal boundary conditions
\begin{equation}\label{nlbcs}
u(0,y,t)=u(1,y,t),\; u_{xx}(0,y,t)=u_{xx}(1,y,t),
\end{equation}
 and initial condition
\begin{equation}\label{ic}
u(x,y,0)=\phi(x,y),\qquad (x,y)\in \Omega:(0,1)\times(0,1),
\end{equation}
where $\Pi:=\Omega \times (0,T),$
${}^cD_{0_+,t}^{\alpha}$ and ${}^cD_{0_+,t}^{\alpha_i}$ stands for left sided Caputo fractional derivative of order $\alpha$ and $\alpha_i$ respectively such that $0<\alpha_m<...<\alpha_1<\alpha\leq 1,\; \; m\in \mathbb{Z}^+,\;\psi_i\geq0 $ and $\Delta^2:=\displaystyle \frac{\partial ^4}{\partial x^4}+\frac{\partial ^4}{\partial y^4}.$

The solution of the inverse problem is usually recovered via an additional information of the system. By additional information, we mean some extra data about the behavior of the system that may be given at intermediate point \cite{ali2time}, on the final time \cite{Kirane} or on the whole domain \cite{ismailov}. Such information/data is termed as over-determination/over-specified condition in the context of inverse problems.

For unique solvability of inverse problem \eqref{problem}-\eqref{ic}, we will consider total energy of the system as over-determination condition, i.e.,
\begin{equation}\label{odc}
\int_{0}^{1}\int_{0}^{1}u(x,y,t)dxdy=E(t).
\end{equation}

History of non integer order integrals and derivatives, termed as fractional calculus (FC), goes back to late seventeenth century.  In the last five decades, FC evolved quite rapidly. Most of the books, journals and conferences were held due to its vast applications in many fields of sciences and engineering are summarized nicely in \cite{History4,History5}.  For detailed evolution of FC, readers are referred \cite{History1,History2,History3}. Differential equations involving the non-integer order derivative of unknown function are termed as Fractional differential equations (FDEs). FDEs have been used extensively to model many well known physical phenomena. Just to mention a few,  FDEs are used to model different processes in viscoelasticity \cite{mainardibook}, image processing
\cite{cuesta}, combustion \cite{Pagnini},  dynamics and control theory \cite{Luo}, motion of Newtonian fluid  \cite{bagley}, earth science dynamics \cite{DLi}, cosmology \cite{Uchaikin}, financial economics \cite{Fallahgoul}, biological processes \cite{Ionescu}.

Much of the work published to date has been concerned with FDEs involving single fractional derivative. For detail, we refer readers to \cite{Diethelm, Kilbas, Miller, Podlubny, Samko} and references cited therein.
 Recently, many authors (see e.g. \cite{Chechkin, Chechkin1, Chechkin2, Naber, Sokolov}) addressed the sub-diffusion processes in which  the mean square displacement of the particles is of logarithmic growth.  To model such ultra slow diffusive processes, scientists have used distributed order fractional derivatives which is an integral of fractional derivatives with respect to continuously changing orders.  Multi-term time-fractional derivative, the special case of distributed order fractional derivative, involves more than one fractional order differential operators. Many real world problems are modelled more appropriately by using multi-term time fractional derivatives.


Let us provide some literature concerning inverse problems for FDEs. Tuan and Nane \cite{Tuan} discussed the inverse problem, of extracting space dependent source term, for time fractional diffusion equation with final temperature over-determination condition. They regularize the solution by using trigonometric method in nonparametric regression. In \cite{Yong}, authors have used revised generalized Tikhonov regularization method to extract the unknown space dependent source term for time fractional diffusion equation. Extraction of space dependent source, by using generalized Fourier's method,  for space time fractional differential equation is considered in \cite{Sarafcaa}. Two inverse problems for fourth-order parabolic fractional differential equation were investigated by Sara et. al in \cite{Sara4thorder}. In \cite{Luchkomulti}, uniqueness and a priory estimates of the solution is established by maximum principle for multi-term time fraction diffusion equation.

Jacobi tau method is used for multi-term time space fraction differential equation with Dirichlet boundary conditions was studied in \cite{Bhrawy}. Analytical and numerical solution for the direct problem for the two dimensional multi-term fractional differential equation with zero Dirichlet boundary conditions were analyzed in \cite{Shujun}. The classical solution of the inverse problem for two dimensional FDE involving a Caputo fractional derivative is determined  in \cite{malik-inv2D}.


In this paper, we will prove the existence of the classical nature of the solution of the inverse problem given by \eqref{problem}-\eqref{odc}. By classical solution of the inverse problem, we mean a pair of function $u(x,y,t)$ and $a(t)$ such that $u(x,y,t) \in C(\Pi),\, \Delta u\in C(\Pi),\, \bigl({}^cD_{0_+,t}^{\alpha}+\sum_{i=1}^{m}\psi_i
{}^cD_{0_+,t}^{\alpha_i}\bigr)\in C(\Pi)$ and $a(t)\in C(0,T).$

The rest of the paper is organized as follows: In Section \ref{Bi-orthogonal System}, bi-orthogonal system consisting of Riesz bases obtained from eigenfunction of the spectral and conjugate problem is presented. Some key lemmata alongside some basic definitions are given in Section \ref{preliminaries}. The research findings of the article are demonstrated in Section \ref{mainresult}.

\section{Bi-orthogonal System}\label{Bi-orthogonal System}
We intend to use Fourier's method for construction of the solution of the inverse problem. Crust of the Fourier's method is to use the eigenfunctions to write the series expansion of the unknown functions. Set of eigenfunctions can only be used if it is complete, minimal and form the basis in some appropriate functional space.

The spectral  and the conjugate problem corresponding to \eqref{problem}-\eqref{nlbcs} are
\begin{equation}\label{sp}
\left.
\begin{array}{lcl}

\vspace{0.3cm}

\displaystyle\frac{\partial ^4 Z}{\partial x^4}(x,y)+\frac{\partial ^4 Z}{\partial x^4}(x,y)=\sigma Z(x,y),\\

\vspace{0.3cm}

\displaystyle\frac{\partial Z}{\partial x}(0,y)=0=\frac{\partial^3 Z}{\partial x^3}(0,y),
\frac{\partial Z}{\partial y}(x,0)=0=\frac{\partial Z}{\partial y}(x,1),\\

\vspace{0.3cm}

\displaystyle\frac{\partial^3 Z}{\partial^3 y}(x,0)=0=\frac{\partial^3 Z}{\partial^3 y}(x,1),\\

\vspace{0.3cm}

Z(0,y)=Z(1,y),\; \displaystyle\frac{\partial^2 Z}{\partial x^2}(0,y)=\frac{\partial^2 Z}{\partial x^2}(1,y),
\end{array}
\right\}
\end{equation}
and
\begin{equation}\label{cp}
\left.
\begin{array}{lcl}

\vspace{0.3cm}

\displaystyle\frac{\partial ^4 W}{\partial x^4}(x,y)+\displaystyle\frac{\partial ^4 W}{\partial x^4}(x,y)=\sigma W(x,y),\\

\vspace{0.3cm}

W(1,y)=0=\displaystyle\frac{\partial^2 W}{\partial x^2}(1,y),
\displaystyle\frac{\partial W}{\partial y}(x,0)=0=\frac{\partial W}{\partial y}(x,1),\\

\vspace{0.3cm}

\displaystyle\frac{\partial^3W}{\partial y^3}(x,0)=0=\frac{\partial^3 W}{\partial y^3}(x,1),\\

\vspace{0.3cm}

\displaystyle\frac{\partial W}{\partial x}(0,y)=\displaystyle\frac{\partial W}{\partial x}(1,y),\; \displaystyle\frac{\partial^3 W}{\partial x^3}(0,y)=\frac{\partial^3 W}{\partial x^3}(1,y),
\end{array}
\right\}
\end{equation}
respectively.

The eigenvalues of the spectral \eqref{sp} and conjugate problem \eqref{cp} are
$$\sigma_{nk}=\mu_k+\lambda_n\quad \mbox{and} \quad  \sigma_{0k}=\mu_k,$$ where $\mu_k=(k\pi)^4$ and $\lambda_n=(2n\pi)^4,\;k,n\in\mathbb{Z}^+.$

The eigenfunctions of spectral and conjugate problem forms the bi-orthogonal system with the following one-one correspondence \cite{AS berdyshev}
\begin{eqnarray}
\biggl\{ \underbrace{Z_{0k}(x,y)}_\downarrow,&  \underbrace{Z_{(2n-1)k}(x,y)}
_\downarrow,\quad \underbrace{Z_{2nk}(x,y)}_\downarrow \biggr\},\label{bios1}\\
\biggl\{\;W_{0k}(x,y),\;&  W_{(2n-1)k}(x,y),\quad W_{2nk}(x,y)\biggr\},\label{bios2}
\end{eqnarray}
where
\begin{gather*}
  Z_{0k}(x,y)=\sqrt{2}\cos(k\pi y),\;Z_{(2n-1)k}(x,y)=\sqrt{2}\cos(2n  \pi x)\cos(k\pi y),\\
  Z_{2nk}(x,y)=\sqrt{2}x\sin(2n\pi x)\cos(k\pi y),
 W_{0k}(x,y)=2\sqrt{2}(1-x)\cos(k\pi y),\\ W_{(2n-1)k}(x,y)=4\sqrt{2}(1-x)\cos(2n  \pi x)\cos(k\pi y),\;\\ W_{2nk}(x,y)=4\sqrt{2}\sin(2n\pi x)\cos(k\pi y).
\end{gather*}
By using the fact that $a^2+b^2\geq2ab,\; \forall \;a,\; b\; \in \mathbb{R}$, we have
\begin{equation}\label{sigmaest}
\frac{1}{\sigma_{nk}}\leq\frac{1}{n^2k^2}\quad\mbox{and}\quad \frac{1}{\mu_{k}}\leq\frac{1}{k^2}, \quad \forall \,n,\,k \in \mathbb{Z}^+.
\end{equation}
\noindent Notice that
\[|Z_{nk}|\leq \sqrt{2},\quad \; |W_{nk}|\leq 4\sqrt{2},\;\quad \forall \;n\in \mathbb{Z}^+ \cup\{0\},\; k\in \mathbb{Z}^+.\]
\begin{lemma}\cite{AS berdyshev}
$Z_{nk}(x,y)$ and $W_{nk}(x,y)$ given by \eqref{bios1} and \eqref{bios2} form Riesz bases of $L_2(\Omega).$
\end{lemma}
\section{Preliminaries}\label{preliminaries}
This section constitutes of some basic definitions and terminologies that will be used through out the article.
\begin{definition}\cite{Kilbas,Samko}
  Let $h\in L^1_{loc}[a,b],\;-\infty < a <z<b < \infty$, then the left sided Riemann-Liouville fractional integral of order $\xi>0$ is defined as
\begin{eqnarray*}
I_{a_+,z}^{\xi}h(z)=\frac{1}{\Gamma(\xi)}
\int_{a}^{z}\frac{h(\tau)}{(z-\tau)^{1-\xi}}
d\tau,\quad z\in (a,b].
\end{eqnarray*}
\end{definition}

\begin{definition}\cite{Kilbas,Samko}
  Let $h\in AC[a,b]$, then the left sided Caputo fractional derivative of order $\xi$; $0<\xi<1$ is defined as
  \begin{eqnarray*}
    ^c{}D_{a_+,z}^\xi h(z)=I_{a_+,z}^{1-\xi}h'(z).
  \end{eqnarray*}
\end{definition}

\begin{lemma}\cite{Ali-Sara1}\label{convolution}
For $g,\,h\; \in C^1[0,b],$ then
\[\frac{d}{dz}(g(z)*h(z))=g(z)h(0)+g(z)*\frac{d}{dz}h(z)=
h(z)g(0)+h(z)*\frac{d}{dz}g(z).\]
\end{lemma}
\begin{definition}\cite{Luchko}
For $\eta,\;\xi_j>0,\; z_j\in \mathbb{C};\; j=1,2,...,n$ multinomial Mittag-Leffler function is defined as
\begin{equation*}E_{(\xi_1,\xi_2,...\xi_n),\eta}(z_1,z_2,...,z_n):=\sum_{k=0}^{\infty}\sum_{\substack{l_1+...+l_n=k\\
                              l_1\geq 0,..., l_n\geq 0}}(k;l_1,...,l_n)\displaystyle\frac{\displaystyle\Pi_{i=1}^{n}z_i^{l_i}}{\Gamma(\eta+
\displaystyle\sum_{i=1}^{n}\xi_il_i)},\end{equation*}
where
$\displaystyle(k;l_1,...,l_n):=\frac{k!}{l_1!\times...\times l_n!}.$
\end{definition}

\begin{remark}
For $z_j=0,\; j=2,3,...,n$ multinomial Mittag-Leffler function be reduced to two parameter Mittag-Leffler function.
\begin{align*}
E_{(\xi_1,\xi_2,...,\xi_n),\eta}(z_1,0,...,0)&=\sum_{k=0}^{\infty}\frac {z_1^k}{\Gamma(\xi_1k+\eta)}:=E_{\xi_1,\eta}(z_1).
\end{align*}
\end{remark}

\begin{remark}\label{symetriceq}
By using definition of multinomial Mittag-Leffler function, we have
\begin{equation*}
 E_{(\xi_1,\xi_1-\xi_2,...,\xi_1-\xi_n),\eta}{(z_1,z_2,...,z_n)}
 =E_{(\xi_1-\xi_n,...,\xi_1-\xi_2,\xi_1),\eta}{(z_n,...,z_2,z_1)}.
\end{equation*}
\end{remark}

\begin{lemma}\label{podlubnyext}
  For $0<\eta<1$ and $0<\xi_n<...<\xi_2<\xi_1<1$ be given. Assume that $\xi_1\pi/2<\mu<\xi_1\pi$, $\mu\leq|arg(z_{n})|\leq\pi$ and $z_i<0, \,i=1,2,...,n$. Then, there exists a constant $C_1$ depending only on $\mu, \xi_i; i=1,2,...,n$ such that
\[|E_{(\xi_1-\xi_n,,...,\xi_1-\xi_2,\xi_1),\eta}(z_n,...,z_2,z_1)|\leq \frac{C_1}{1+|z_{n}|}\leq C_1.\]
\begin{proof} From Remark \ref{symetriceq}, we have
 $$E_{(\xi_1-\xi_n,,...,\xi_1-\xi_2,\xi_1),\eta}{(z_n,...,z_2,z_1)}
 =E_{(\xi_1,\xi_1-\xi_2,...,\xi_1-\xi_n),\eta}{(z_1,z_2,...,z_n)}.$$
Due to Lemma 3.2 of \cite{YAMAMOTOLemma} , we can obtain the required result.
\end{proof}
\end{lemma}
\bigskip

\noindent For convenience, we introduce the following notation
\begin{align*}&\mbox{\Large $e$ \hspace{-.3cm} }_{(m_1\xi_1,m_2\xi_2,...,m_n\xi_n),\eta}(z):=\\&
\hspace{1.5cm}z^{\eta-1}E_{(\xi_1,\xi_2,...,\xi_n),\eta}
(-m_1z^{\xi_1},-m_2z^{\xi_2},-m_3z^{\xi_3},...,-m_nz^{\xi_n}).\end{align*}

\begin{lemma}\label{etaplus1}
For $\eta,\;\xi_j>0,\; m_j>0\,; \, j\in \mathbb{Z}^+,$ we have
\[\int_{0}^{z}\mbox{\Large $e$ \hspace{-.3cm} }_{(m_1\xi_1,m_2\xi_2,...,m_n\xi_n),\eta}(\tau)\;d\tau=
\mbox{\Large $e$ \hspace{-.3cm} }_{(m_1\xi_1,m_2\xi_2,...,m_n\xi_n),\eta+1}(z).\]
\end{lemma}
\begin{proof}
By using defintion of $\mbox{\Large $e$ \hspace{-.3cm} }_{(m_1\xi_1,m_2\xi_2,...,m_n\xi_n),\eta:\tau}$ and term by term integration, required result can be obtained.
\end{proof}
\begin{remark}
For $m_j=0,\;j=2,...,n$, in Lemma \ref{etaplus1}, we get \begin{eqnarray*}
\int_{0}^{z}\mbox{\Large $e$ \hspace{-.3cm} }_{(m_1\xi_1,0,0...,0),\eta}(\tau)d\tau=
\int_{0}^{z}\tau^{\eta-1}{E}_{\xi_n,\eta}(-m_n \tau)d\tau=z^\eta{E}_{\xi_1,\eta+1}(-m_1z^{\xi_1}).
\end{eqnarray*}
Above relation was mentioned in \cite{Sara4thorder}.
\end{remark}

\begin{lemma}\label{conv}
  For $h\in C^{1}[0,b]$ and $\xi_i>0$, $\psi_i>0$, for $i=1,2,...,n$, we have
$$|h(z)\ast\mbox{\Large $e$ \hspace{-.2cm}} _{(m_n(\xi_1-\xi_n),...,m_2(\xi_1-\xi_2),m_1\xi_1),\xi_1}(z)|
\leq\frac{C_{1}}{m_{n}}\| h\|_1,$$
where $\|h\|_1=\displaystyle\max_{a<z<b}|h(z)|$ and $``\ast"$ represents the Laplace convolution.
\begin{proof}
By using the fact that $h\in C^{1}[a,b]$ alongside Lemma \ref{etaplus1} and Lemma \ref{podlubnyext}, we can obtain required result.
\end{proof}
\end{lemma}
\begin{remark}
In \cite{Ali-Sara1}, upper bound for the convolution of the Mittag-Leffler function with any continuously differentiable function is established, which is given by
\begin{equation}\label{conold}
|g(z)*z^{\xi_1}E_{\xi_1,\xi_1+1}(-m_nz^{\xi_1})|\leq \frac{C_1}{m_n}\|g(z)\|.
\end{equation}
Equation \eqref{conold} can also be established by setting $m_{i}=0,\,i=1,...,n-1$ where $n \in N$ in Lemma.

\end{remark}

Some formulae for Laplace transform are mentioned here as we intend to use
Laplace transform to solve system of multi-term time fractional differential equations.

\begin{itemize}
\item $
\displaystyle \mathcal{L}\bigl({}^cD_{0_+,z}^{\eta}+\sum_{i=1}^{m}\psi_i
{}^cD_{0_+,z}^{\eta_i}\bigr)h(z)=(s^{\eta}+\sum_{i=1}^{m}
\psi_is^{\eta_i})
\mathcal{L}(h(z))-(s^{\eta-1}+\sum_{i=1}^{m}\psi_is^{\eta_i-1})
h(0).
$
\item $
\displaystyle\mathcal{L}\bigl(\mbox{\Large $e$ \hspace{-.3cm} }_{(m_1\xi_1,m_2\xi_2,...,m_n\xi_n),\eta}(z)\bigr)
=\frac{s^{-\eta}}{1+\displaystyle\sum_{i=2}^{n}m_is^{-\xi_i}}.
$
\end{itemize}
\section{Inverse Source Problem}\label{mainresult}
In this section, first we are going to present some relevant lemmata that would be advantageous in ultimately proving the main result i.e., the classical nature of the solution.
Certain regularity conditions are imposed on given datum to prove the well-posedness of the inverse problem \eqref{problem}-\eqref{odc}.

\begin{lemma}\label{fest}
For $h\in C^{2,1,0}_{x,y,t}(\Pi)$ such that $h(0,y,t)=h(1,y,t),$
then we have
\begin{itemize}
\item $|h_{0k}(t)|\leq \displaystyle\frac{2}{k \pi}\bigg|\left\langle \frac{\partial h}{\partial y},\sqrt{2}\sin(k \pi y)\right\rangle\bigg|,$
\item $|h_{(2n-1)k}(t)|\leq \displaystyle \frac{4}{(2n \pi)(k \pi)}
\bigg|\left\langle \frac{\partial^2h}{\partial x \partial y},\sqrt{2}\sin(2n \pi)\sin(k \pi y)\right\rangle\bigg|,$
\item $|h_{2nk}(t)|\leq \displaystyle \frac{4}{(2n \pi)(k \pi)}
\bigg|\left\langle \frac{\partial^2 h}{\partial x\partial y},\sqrt{2}\cos(2n \pi)\sin(k \pi y)\right\rangle\bigg|$,
\item $|h_{(2n-1)k}(t)|\leq \displaystyle \frac{4}{(2n \pi)^2}\bigg\{\bigg|\left\langle\frac{\partial h(1,y)}{\partial x}+\frac{\partial h(0,y)}{\partial x},\sqrt{2}\cos(k \pi y)\right\rangle\bigg|+\bigg|\left\langle \frac{\partial^2h}{\partial x^2 },\sqrt{2}\cos(2n \pi)\cos(k \pi y)\right\rangle\bigg|\bigg\},$
\item $|h_{2nk}(t)|\leq \displaystyle\frac{4}{(2n \pi)^2}\bigg|\left\langle \frac{\partial^2 h}{\partial x^2},\sqrt{2}\sin(2n \pi)\cos(k \pi y)\right\rangle\bigg|,$
 \end{itemize}
 where $h_{nk}(t)=\left\langle h(x,y,t), W_{nk}(x,y)\right\rangle,\quad n\in \mathbb{Z}^+\cup \{0\},\; k\in \mathbb{Z}^+.$
 \begin{proof}
We are going to prove the fourth inequality,
other inequalities can be obtained in a similar way.

By definition
\[|h_{(2n-1)k}|=|\left\langle h(x,y), W_{(2n-1)k}\right\rangle|\leq 4\sqrt{2}\bigg|\int_{0}^{1}\int_{0}^{1}h(x,y)\cos(2n  \pi x)\cos(k \pi y)dxdy\bigg|,\]
where we have used the fact that $|1-x|<1$ in $\Omega.$

Integration by parts leads to
\[|h_{(2n-1)k}|\leq \frac{4\sqrt{2}}{2n \pi}\bigg|\int_{0}^{1}\int_{0}^{1}h_x(x,y)\sin(2n  \pi x)\cos(k \pi y)dxdy\bigg|.\]
Required inequality can be obtained by performing integration by parts with respect to $x$, once more.
\end{proof}

\end{lemma}

\begin{lemma}\label{phiest}
For $g\in C^{5,5}_{x,y}(\Omega),\, i=0,1,2,3$
and $ j=1,3$ such that \[ \frac{\partial^i g}{\partial x^i}(0,y)= \frac{\partial^i g}{\partial x^i}(1,y),\quad \frac{\partial^j g}{\partial x^j}(x,0)=0= \frac{\partial^j g}{\partial x^j}(x,1), \]
then, we have
\begin{itemize}
\item $|g_{0k}|\leq \displaystyle\frac{2}{(k \pi)}\bigg|\left\langle \frac{\partial g}{\partial y},\sqrt{2}\sin(k \pi y)\right\rangle\bigg|,$
\item $|g_{0k}|\leq \displaystyle\frac{2}{(k \pi)^5}\bigg|\left\langle \frac{\partial^5 g}{\partial y^5},\sqrt{2}\sin(k \pi y)\right\rangle\bigg|,$
\item $|g_{(2n-1)k}|\leq \displaystyle \frac{4}{(2n \pi )(k \pi)}
\bigg|\left\langle \frac{\partial^2g}{\partial x \partial y},\sqrt{2}\sin(2n \pi )\sin(k \pi y)\right\rangle\bigg|,$
\item $|g_{(2n-1)k}|\leq \displaystyle \frac{4}{(2n \pi )^2}\bigg\{\bigg|\left\langle\frac{\partial g(1,y)}{\partial x}+\frac{\partial g(0,y)}{\partial x},\sqrt{2}\cos(k \pi y)\right\rangle\bigg|
+\bigg|\left\langle \frac{\partial^2g}{\partial x^2 },\sqrt{2}\cos(2n \pi x)\cos(k \pi y)\right\rangle\bigg|\bigg\},$
\item $|g_{(2n-1)k}|\leq \displaystyle \frac{4}{(2n \pi )^5(k \pi)}\left\{\bigg|\left\langle\frac{\partial^6 g}{\partial^5 x \partial y},\sqrt{2}\sin(2n \pi x)\sin(k \pi y)\right\rangle\bigg|\right\},$
\item $|g_{(2n-1)k}|\leq \displaystyle \frac{4}{(2n \pi )(k \pi)^5}\left\{\bigg|\left\langle\frac{\partial^6 g}{\partial^5 y \partial x},\sqrt{2}\sin(2n \pi x)\sin(k \pi y)\right\rangle\bigg|\right\},$
\item $|g_{2nk}|\leq \displaystyle \frac{4}{(2n \pi )(k \pi)}
\bigg|\left\langle \frac{\partial^2 g}{\partial x\partial y},\sqrt{2}\cos(2n \pi x)\sin(k \pi y)\right\rangle\bigg|$,
\item $|g_{2nk}|\leq \displaystyle\frac{4}{(2n \pi )^2}\bigg|\left\langle \frac{\partial^2 g}{\partial x^2},\sqrt{2}\sin(2n \pi x)\cos(k \pi y)\right\rangle\bigg|,$
\item $|g_{2nk}|\leq \displaystyle\frac{4}{(2n \pi )^5(k \pi)}\bigg\{\bigg|\left\langle \frac{\partial^6 g(1,y)}{\partial x^5\partial y}- \frac{\partial^6 g(0,y)}{\partial x^5\partial y},\sqrt{2}\sin(k \pi y)\right\rangle\bigg|+\bigg|\left\langle \frac{\partial^6 g}{\partial x^5\partial y},\sqrt{2}\cos(2n \pi x)\sin(k \pi y)\right\rangle\bigg|\bigg\},$
\item $|g_{2nk}|\leq \displaystyle\frac{4}{(2n \pi )(k \pi)^5}\bigg\{\bigg|\left\langle \frac{\partial^5 g(1,y)}{\partial y^5}- \frac{\partial^5 g(0,y)}{\partial y^5},\sqrt{2}\sin(k \pi y)\right\rangle\bigg|+\bigg|\left\langle \frac{\partial^6 g}{\partial x\partial y^5},\sqrt{2}\cos(2n \pi x)\sin(k \pi y)\right\rangle\bigg|\bigg\},$
 \end{itemize}
where $g_{nk}=\left\langle g(x,y), W_{nk}\right\rangle,\quad n\in \mathbb{Z}^+\cup \{0\},\; k\in \mathbb{Z}^+.$
\end{lemma}

\begin{theorem}
 Under the following regularity conditions on the given datum $\phi(x,y),\; f(x,y,t)$ and $E(t)$ the solution of inverse source problem \eqref{problem}-\eqref{odc} is classical in nature.
  \begin{itemize}
  \item $\phi(x,y)\in C^{5,5}_{x,y}(\Omega),\,i=0,1,2,3,$ and  $j=1,3$ such that
\[ \frac{\partial^i \phi}{\partial x^i}(0,y)= \frac{\partial^i \phi}{\partial x^i}(1,y),\quad \frac{\partial^j \phi}{\partial x^j}(x,0)=0= \frac{\partial^j\phi}{\partial x^j}(x,1),\]
\item $f(x,y,t) \in C^{2,1,0}_{x,y,t}(\Pi)$ such that $f(0,y,t)=f(1,y,t),$
 \item  $\biggl(\displaystyle\int_{0}^{1}\int_{0}^{1}f(x,y,t)dxdy\biggr)^{-1}\leq C_2,$ for some $C_2>0,$
\item $E(t)\in AC[0,T]$ and $\displaystyle\int_{0}^{1}\int_{0}^{1}\phi(x,y)dxdy=E(0).$
\end{itemize}
\end{theorem}
\begin{proof}
In order to prove the classical nature of the solution of inverse source problem \eqref{problem}-\eqref{odc}, firstly we will construct the solution of the problem followed by existence, uniqueness and stability results.

\textbf{Construction of the solution:}\\
We seek the series representation of the solution by making use of bi-orthogonal system \eqref{bios1}-\eqref{bios2} to get
\begin{align}
  u(x,y,t)&=\sum_{k=1}^{\infty}T_{0k}(t)Z_{0k}(x,y)+
  \sum_{n=1}^{\infty}
  \sum_{k=1}^{\infty}\big[T_{(2n-1)k}(t)Z_{(2n-1)k}(x,y)+T_{2nk}(t)Z_{2nk}(x,y)\big],\label{seriesu}
\end{align}
where $T_{nk}$'s satisfy the following multi-term time fractional differential system
\begin{gather}
 \bigl({}^cD_{0_+,t}^{\alpha}+\sum_{i=1}^{m}\psi_i
{}^cD_{0_+,t}^{\alpha_i}\bigr)T_{0k}(t)=-\mu_kT_{0k}(t)+a(t)f_{0k}(t), \label{fdesT0} \\
 \bigl({}^cD_{0_+,t}^{\alpha}+\sum_{i=1}^{m}\psi_i
{}^cD_{0_+,t}^{\alpha_i}\bigr)T_{(2n-1)k}(t)
  =-\sigma_{nk}T_{(2n-1)k}(t)+4(\lambda_n)^{3/4}
  T_{2nk}(t)+a(t)f_{(2n-1)k}(t),\\
 \bigl({}^cD_{0_+,t}^{\alpha}+\sum_{i=1}^{m}\psi_i
{}^cD_{0_+,t}^{\alpha_i}\bigr)T_{2nk}(t) =-\sigma_{nk}T_{2nk}(t)+a(t)f_{2nk}(t)\label{fdesT2nk},
  \end{gather}
where $f_{nk}(t)=\left\langle f(x,t), W_{nk}(x,y)\right\rangle,\quad n\in \mathbb{Z}^+\cup\{0\},\; k\in \mathbb{Z}^+.$

The solution of the equations \eqref{fdesT0}-\eqref{fdesT2nk} is obtained by using Laplace transform and is given by
\begin{align}
T_{0k}(t)=&\mbox{\Large $e$ \hspace{-.3cm} }_{(\psi_1(\alpha-\alpha_1),...,\psi_m(\alpha-\alpha_m),
\mu_k\alpha),1}(t)\;\phi_{0k}+
\sum_{i=1}^{m}\psi_i\mbox{\Large $e$ \hspace{-.3cm} }_{(\psi_1(\alpha-\alpha_1),...,
\psi_m(\alpha-\alpha_m),\mu_k\alpha_i),\alpha+1-\alpha_1}(t)\;\phi_{0k}\nonumber\\&+a(t)f_{0k}(t)* \mbox{\Large $e$ \hspace{-.3cm} }_{(\psi_1(\alpha-\alpha_1),...,\psi_m(\alpha-\alpha_m),\mu_k\alpha_i),
\alpha}(t),\label{T0k}
\end{align}
\begin{align}
T_{2nk}(t)=&\mbox{\Large $e$ \hspace{-.3cm} }_{(\psi_1(\alpha-\alpha_1),...,\psi_m(\alpha-\alpha_m),\sigma_{nk}
\alpha),1}(t)\;\phi_{2nk}
+\sum_{i=1}^{m}\psi_i\mbox{\Large $e$ \hspace{-.3cm} }_{(\psi_1(\alpha-\alpha_1),...,
\psi_m(\alpha-\alpha_m),\sigma_{nk}\alpha_i),\alpha+1-\alpha_i}(t)\;\phi_{2nk}
\nonumber\\&+a(t)f_{2nk}(t)*
\mbox{\Large $e$ \hspace{-.3cm} }_{(\psi_1(\alpha-\alpha_1),...,\psi_m(\alpha-\alpha_m),
\sigma_{nk}\alpha_i),\alpha}(t),\label{T2nk}\end{align}
\begin{align}
T_{(2n-1)k}(t)=&\mbox{\Large $e$ \hspace{-.3cm} }_{(\psi_1(\alpha-\alpha_1),...,\psi_m(\alpha-\alpha_m),\sigma_{nk}
\alpha),1}(t)\;\phi_{(2n-1)k}+a(t)f_{(2n-1)k}(t)
*\mbox{\Large $e$ \hspace{-.3cm} }_{(\psi_1(\alpha-\alpha_1),...,\psi_m(\alpha-\alpha_m),\sigma_{nk}
\alpha_i),\alpha}(t)\nonumber\\&+\sum_{i=1}^{m}\psi_i
\mbox{\Large $e$ \hspace{-.3cm} }_{(\psi_1(\alpha-\alpha_1),...,
\psi_m(\alpha-\alpha_m),\sigma_{nk}\alpha_i),\alpha+1-\alpha_i}(t)\;\phi_{(2n-1)k}\nonumber\\&
+4(\lambda_n)^{3/4}T_{2nk}(t)*
\mbox{\Large $e$ \hspace{-.3cm} }_{(\psi_1(\alpha-\alpha_1),...,\psi_m(\alpha-\alpha_m),
\sigma_{nk}\alpha_i),\alpha}(t)\label{T2n1k},
\end{align}
where $a(t)$ is still to be determined.

After getting the series expression of $u(x,y,t)$, we will discuss the existence of the source term $a(t)$.
Use of the over-determination condition \eqref{odc}, for determination of $a(t)$,
\begin{eqnarray}
\int_{0}^{1}\int_{0}^{1}\bigl({}^cD_{0_+,t}^{\alpha}+\sum_{i=1}^{m}\psi_i
{}^cD_{0_+,t}^{\alpha_i}\bigr)u(x,y,t)dxdy
=\bigl({}^cD_{0_+,t}^{\alpha}+\sum_{i=1}^{m}\psi_i
{}^cD_{0_+,t}^{\alpha_i}\bigr)E(t).
\end{eqnarray}
By using \eqref{problem}, we have
\begin{eqnarray}
\int_{0}^{1}\int_{0}^{1}\bigl(-\frac{\partial ^4u}{\partial x^4}-\frac{\partial ^4u}{\partial y^4}+a(t)f(x,y,t)\bigr)dxdy
=\bigl({}^cD_{0_+,t}^{\alpha}+\sum_{i=1}^{m}\psi_i
{}^cD_{0_+,t}^{\alpha_i}\bigr)E(t).
\end{eqnarray}
Consequently, we have following explicit expression of $a(t)$
\begin{equation}\label{a(t)}
a(t)=\left(\int_{0}^{1}\int_{0}^{1}f(x,y,t)dxdy\right)^{-1}
\bigl({}^cD_{0_+,t}^{\alpha}+\sum_{i=1}^{m}\psi_i
{}^cD_{0_+,t}^{\alpha_i}\bigr)E(t).
\end{equation}
The solution of the inverse source problem is given by \eqref{seriesu} and \eqref{a(t)} where $T_{0k}(t),\; T_{(2n-1)k}(t)$ and $T_{2nk}(t)$ are given by \eqref{T0k}-\eqref{T2n1k}.

 Next, we will show that the obtained solution is regular in nature.

\noindent\textbf{Existence of the solution:}\\
We will present the existence result for the solution of the inverse source problem i.e., $\{a(t),u(x,y,t)\}.$

\noindent Under the given assumptions on $f(x,y,t)$ and $E(t), a(t)\in C[0,T].$

For the existence of the solution \eqref{seriesu}, we need to show that the series representation of $u(x,y,t)$, $\Delta^2u(x,y,t)$ and $\bigl({}^cD_{0_+,t}^{\alpha}+\displaystyle\sum_{i=1}^{m}\psi_i
{}^cD_{0_+,t}^{\alpha_i}\bigr)u(x,y,t)$ are uniformly convergent.

\begin{itemize}
\item \textbf{Uniform convergence of series representation of $u(x,y,t)$}
\end{itemize}
In order to establish the uniform convergence of series representation of $u(x,y,t)$ we will show that the double infinite series $\displaystyle\sum_{k=1}^{\infty} T_{0k}(t),\;\displaystyle\sum_{n=1}^{\infty}\sum_{k=1}^{\infty} T_{nk}(t)$ are uniformly convergent.

From equation \eqref{T0k}, we get
\begin{align}
|T_{0k}(t)|\leq&|\mbox{\Large $e$ \hspace{-.3cm} }_{(\psi_1(\alpha-\alpha_1),...,\psi_m(\alpha-\alpha_m),
\mu_k\alpha),1}(t)\;\phi_{0k}|\nonumber\\&+
\sum_{i=1}^{m}\psi_i|\mbox{\Large $e$ \hspace{-.3cm} }_{(\psi_1(\alpha-\alpha_1),...,
\psi_m(\alpha-\alpha_m),\mu_k\alpha_i),\alpha+1-\alpha_i}(t)\;\phi_{0k}|\nonumber\\&+|a(t)f_{0k}* \mbox{\Large $e$ \hspace{-.3cm} }_{(\psi_1(\alpha-\alpha_1),...,\psi_m(\alpha-\alpha_m),\mu_k\alpha_i),\alpha}(t)|.\label{T0kest1}
\end{align}
Since, $a(t)\in C[0,T]$ so, $\exists$ a constant $C_3$ such that  $|a(t)|\leq C_3,\; \forall \;t\in (0,T].$ By using the fact that $|a(t)|\leq C_3$. By using Lemma \ref{podlubnyext} and Lemma \ref{conv}, we get
\begin{align*}
|T_{0k}(t)|\leq {C_1}|\phi_{0k}|(1+\sum_{i=1}^{m}\psi_i\ T^{\alpha-\alpha_i})+\frac{C_1C_3}{(k \pi)^4}|f_{0k}(t)|.
\end{align*}
By using Lemma \ref{phiest}, we have
\[|\phi_{0k}|\leq \frac{2}{k \pi}\bigg|\left\langle\frac{\partial \phi}{\partial y},\sqrt{2}\sin(k \pi y)\right\rangle\bigg|.\]
Since, $ab\leq \frac{1}{2}(a^2+b^2),\; \forall \; a,b\in \mathbb{R},$ we have
\[|\phi_{0k}|\leq \frac{1}{(k \pi)^2}+\bigg|\left\langle\frac{\partial \phi}{\partial y},\sqrt{2}\sin(k \pi y)\right\rangle\bigg|^2.\]
From Equation \eqref{T0kest1}, we have
\begin{align*}
\bigg|\sum_{k=1}^{\infty}T_{0k}(t)\bigg|\leq& {C_1}(1+\sum_{i=1}^{m}\psi_i\ T^{\alpha-\alpha_i})
\sum_{k=1}^{\infty}\bigg(\frac{1}{(k \pi)^2}+\big|\left\langle\frac{\partial \phi}{\partial y},\sqrt{2}\sin(k \pi y)\right\rangle\big|^2\bigg)\nonumber\\&+\sum_{k=1}^{\infty}\frac{C_1C_3}{(k \pi)^4}|f_{0k}(t)|.
\end{align*}
As $\displaystyle\{\sqrt{2}\sin(k \pi y)\}_{k=1}^{\infty}$ forms an orthonormal sequence in $L_2(0,1)$, by Bessel Inequality, we have
\begin{align*}
\bigg|\sum_{k=1}^{\infty}T_{0k}(t)\bigg|\leq {C_1}(1+\sum_{i=1}^{m}\psi_i\ T^{\alpha-\alpha_i})
\bigg(\sum_{k=1}^{\infty} \frac{1}{(k \pi)^2}+\bigg\|\frac{\partial \phi}{\partial y}\bigg\|_2\bigg)+\sum_{k=1}^{\infty}\frac{C_1C_3}{(k \pi)^4}\|f\|_3,
\end{align*}
where $\|\phi\|_2=\sqrt{\left\langle \phi(x,y),\,\phi(x,y)\right\rangle}$ and $\|f\|_3=\sqrt{\left\langle f(x,y,t),\,f(x,y,t)\right\rangle}$.

By using  Cauchy-Bunyakovsky-Schwarz Inequality and the fact that $|W_{0k}|\leq 4\sqrt{2}$, we get
\begin{align}\label{T0kest}
\bigg|\sum_{k=1}^{\infty}T_{0k}(t)\bigg|\leq {C_1}(1+\sum_{i=1}^{m}\psi_i\ T^{\alpha-\alpha_i})
\bigg(\sum_{k=1}^{\infty} \frac{1}{(k \pi)^2}+\bigg\|\frac{\partial \phi}{\partial y}\bigg\|_2\bigg)+\sum_{k=1}^{\infty}\frac{4\sqrt{2}C_1C_3}{(k \pi)^4}\|f\|_3.
\end{align}


\noindent Similarly, for $T_{2nk}(t),$ we have
\begin{align}\label{T2nkest}
\bigg|\sum_{n=1}^{\infty}\sum_{k=1}^{\infty}T_{2nk}(t)\bigg|&\leq 2{C_1}(1+\sum_{i=1}^{m}\psi_i\ T^{\alpha-\alpha_i})
\bigg(\sum_{n=1}^{\infty}\sum_{k=1}^{\infty} \frac{1}{(2n \pi)^2(k \pi)^2}+\bigg\|\frac{\partial^2 \phi}{\partial x\partial y}\bigg\|_2\bigg)\nonumber\\&+\sum_{n=1}^{\infty}\sum_{k=1}^{\infty}\frac{2\sqrt{2}C_1C_3}{(2n \pi)^2(k \pi)^2}\|f\|_3.
\end{align}

For uniform convergence of series involved in $u(x,y,t)$ we also need to analyze $T_{(2n-1)k}(t).$
On same lines as for $T_{0k}(t)$ and $T_{2nk}(t)$, \eqref{T2n1k} takes the form
\begin{align*}
\bigg|\sum_{n=1}^{\infty}\sum_{k=1}^{\infty}T_{(2n-1)k}(t)\bigg|
\leq &\sum_{n=1}^{\infty}\sum_{k=1}^{\infty}\bigg[\frac{4{C_1}(1+\sum_{i=1}^{m}\psi_i\ T^{\alpha-\alpha_i})}{(2n \pi)(k \pi)}
\langle \frac{\partial^2 \phi}{\partial x \partial y}, \sqrt{2}\sin(2n  \pi x)\sin(k  \pi y)\rangle\\&+\frac{C_1C_3}{2(2n \pi)^2(k \pi)^2}|f_{(2n-1)k}|\nonumber+\frac{8{C_1}(1+\sum_{i=1}^{m}\psi_i\ T^{\alpha-\alpha_i})}{(2n \pi)(k \pi)^2}
\langle \frac{\partial^2 \phi}{\partial x^2}, \sqrt{2}\sin(2n  \pi x)\cos(k  \pi y)\rangle\\&+\frac{C_1C_3}{2(2n \pi)(k \pi)^4}|f_{2nk}|\bigg].
\end{align*}
Since, $\displaystyle\{\sqrt{2}\sin(2n  \pi x)\sin(k  \pi y)\}_{n,k=1}^{\infty},$ and $\displaystyle\{\sqrt{2}\sin(2n  \pi x)\cos(k  \pi y)\}_{n,k=1}^{\infty}$ form an orthonormal sequences in $L_2(\Omega)$,  by using the elementary inequality $2ab<a^2+b^2$, the Bessel
and the Cauchy-Bunyakovsky-Schwarz Inequality, we have
\begin{align}\label{T2n1kest}
\bigg|\sum_{n=1}^{\infty}\sum_{k=1}^{\infty}T_{(2n-1)k}(t)\bigg|
\leq &2{C_1}(1+\sum_{i=1}^{m}\psi_i\ T^{\alpha-\alpha_i})
\bigg(\sum_{n=1}^{\infty}\sum_{k=1}^{\infty} \frac{1}{(2n \pi)^2(k \pi)^2}\nonumber+\bigg\|\frac{\partial^2 \phi}{\partial x\partial y}\bigg\|_2\bigg)\\&+\sum_{n=1}^{\infty}\sum_{k=1}^{\infty}\frac{2\sqrt{2}C_1C_3}{(2n \pi)^2(k \pi)^2}\|f\|_3\nonumber\\&+2{C_1}(1+\sum_{i=1}^{m}\psi_i\ T^{\alpha-\alpha_i})
\bigg(\sum_{n=1}^{\infty}\sum_{k=1}^{\infty} \frac{1}{(2n \pi)^2(k \pi)^4}\nonumber+\bigg\|\frac{\partial^2 \phi}{\partial x^2}\bigg\|_2\bigg)\\&+C_1C_3\bigg(\sum_{n=1}^{\infty}\sum_{k=1}^{\infty} \frac{1}{(2n \pi)^2(k \pi)^8}+\big\|f\big\|^2_2\bigg).
\end{align}


Inequalities \eqref{T0kest}-\eqref{T2n1kest} along side  Weierstrass M test ensure that series representation of $u(x,y,t)$ is uniformly convergent. Hence, $u(x,y,t)$ represents a continuous function.

\begin{itemize}
\item \textbf{Uniform convergence of series representation of $\Delta^2u(x,y,t)$}
\end{itemize}
\smallskip

\noindent Note that
\begin{align*}\Delta^2u(x,y,t)&=\sum_{k=1}^{\infty}\mu_kT_{0k}(t)Z_{0k}(x,y)+
  \sum_{n=1}^{\infty}
  \sum_{k=1}^{\infty}\big[\sigma_{nk}T_{(2n-1)k}(t)Z_{(2n-1)k}(x,y)
  \\&+\sigma_{nk}T_{2nk}(t)Z_{2nk}(x,y)\big].\end{align*}
  Since $|Z_{nk}|\leq \sqrt{2},\; \forall n\in \mathbb
\mathbb{Z}^+\cup\{0\},\; k\in \mathbb{Z}^+$, so in order to establish the convergence of $\Delta^2u(x,y,t),$ we need to prove the convergence of
  \[\sum_{k=1}^{\infty}\mu_kT_{0k}(t),\qquad \sum_{n=1}^{\infty}
  \sum_{k=1}^{\infty}\sigma_{nk}T_{(2n-1)k}(t),\qquad \sum_{n=1}^{\infty} \sum_{k=1}^{\infty}
  \sigma_{nk}T_{2nk}(t).\]
By using \eqref{T0k} and Lemma \ref{conv}, we have
\[\big|\sum_{k=1}^{\infty}\mu_kT_{0k}(t)\big|\leq\sum_{k=1}^{\infty} {C_1}\mu_k(1+\sum_{i=1}^{m}\psi_i\ T^{\alpha-\alpha_i})|\phi_{0k}|+C_1C_3 |f_{0k}(t)| .\]
By using Lemma \ref{fest}, Lemma \ref{phiest} and Bessel Inequality, we get
\begin{align*}
  \sum_{k=1}^{\infty}\mu_{k}|T_{0k}(t)|\leq &C_1(1+\sum_{i=1}^{m}\psi_i\ T^{\alpha-\alpha_i})\bigg(\sum_{k=1}^{\infty}\frac{1}{(k \pi)^2}+\bigg\|\frac{\partial^5\phi}{\partial y^5}\bigg\|_2^2\bigg)
\\&+C_1C_3\bigg(\sum_{k=1}^{\infty}\frac{1}{(k \pi)^2}+\big\|\frac{\partial f}{\partial y}\big\|_3^2\bigg).
\end{align*}
We will now turn our attention to $\sigma_{nk}T_{2nk}(t)$, using \eqref{T2nk} to have
\begin{align*}
 \sum_{n=1}^{\infty}
  \sum_{k=1}^{\infty}\sigma_{nk}|T_{2nk}(t)|&\leq\sum_{k=1}^{\infty} {C_1}(1+\sum_{i=1}^{m}\psi_i\ T^{\alpha-\alpha_i})((2n \pi)^4+(k \pi)^4)|\phi_{2nk}|\\&+C_1C_3 |f_{2nk}(t)|,
\end{align*}
due to Lemma  \ref{fest} and \ref{phiest}, we have
\begin{align*}
 \sum_{n=1}^{\infty}
  \sum_{k=1}^{\infty}\sigma_{nk}|T_{2nk}(t)|\leq&\sum_{n=1}^{\infty}\sum_{k=1}^{\infty}\frac{4C_1(1+\sum_{i=1}^{m}\psi_i\ T^{\alpha-\alpha_i})}{(2n \pi)(k \pi)}\bigg\{\bigg|\left\langle \frac{\partial^6 \phi(1,y)}{\partial x^5\partial y}- \frac{\partial^6 \phi(0,y)}{\partial x^5\partial y},\sqrt{2}\sin(k \pi y)\right\rangle\bigg|\\&+\bigg|\left\langle \frac{\partial^6 g}{\partial x^5\partial y},\sqrt{2}\cos(2n \pi x)\sin(k \pi y)\right\rangle\bigg|\bigg\}\\&+\sum_{n=1}^{\infty}\sum_{k=1}^{\infty}\frac{4C_1(1+\sum_{i=1}^{m}\psi_i\ T^{\alpha-\alpha_i})}{(2n \pi)(k \pi)}\bigg\{\bigg|\left\langle \frac{\partial^5 \phi(1,y)}{\partial y^5}- \frac{\partial^5 \phi(0,y)}{\partial y^5},\sqrt{2}\sin(k \pi y)\right\rangle\bigg|\\&+\bigg|\left\langle \frac{\partial^6 \phi}{\partial x\partial y^5},\sqrt{2}\cos(2n \pi x)\sin(k \pi y)\right\rangle\bigg|\bigg\}\\&+\sum_{n=1}^{\infty}\sum_{k=1}^{\infty}\frac{4C_1C_3}{(2n \pi)(k \pi)}\left\langle f_{xy},\sqrt{2}\cos(2n  \pi x)\sin(k \pi y)\right\rangle,
\end{align*}
using the inequality $2ab\leq a^2+b^2$ and the fact  that $\displaystyle\{\sqrt{2}\sin(k \pi y)\}_{k=1}^{\infty}$  and $\displaystyle\{\sqrt{2}\cos(2n \pi x) \sin(k \pi y)\}_{n, k=1}^{\infty}$ form orthonormal sequences in $L_2(0,1)$ and $L_2(\Omega)$ respectively, together with Bessel's Inequality,  we get
\begin{align*}
 \sum_{n=1}^{\infty} \sum_{k=1}^{\infty}\sigma_{nk}|T_{2nk}(t)|\leq&4C_1(1+\sum_{i=1}^{m}\psi_i\ T^{\alpha-\alpha_i})\bigg\{\sum_{n=1}^{\infty}\sum_{k=1}^{\infty}\frac{1}{(2n \pi)^2(k \pi)^2}\\&+\bigg\|\frac{\partial^6 \phi(1,y)}{\partial x^5\partial y}- \frac{\partial^6 \phi(0,y)}{\partial x^5\partial y}\bigg\|_2^2\\&+\bigg\| \frac{\partial^6 \phi}{\partial x^5\partial y}\bigg\|_2^2+\bigg\| \frac{\partial^5 \phi(1,y)}{\partial y^5}- \frac{\partial^5 \phi(0,y)}{\partial y^5}\bigg\|_2^2+\bigg\|\frac{\partial^6 \phi}{\partial x\partial y^5}\bigg\|_2^2\bigg\}\\&+2C_1C_3\sum_{n=1}^{\infty}\sum_{k=1}^{\infty}\bigg\{\frac{1}{(2n \pi)^2(k \pi)^2}+\bigg\| \frac{\partial^2f}{\partial x \partial y}\bigg\|_3^2\bigg\}.
\end{align*}

For the convergence  of $\sigma_{nk}T_{(2n-1)k}(t)$,  consider \eqref{T2n1k} together with Lemma \ref{phiest}, we have
\begin{align*}
 \sum_{n=1}^{\infty}\sum_{k=1}^{\infty}\sigma_{nk}|T_{(2n-1)k}(t)|&\leq 4{C_1}(1+\sum_{i=1}^{m}\psi_i\ T^{\alpha-\alpha_i})\sum_{n=1}^{\infty}\sum_{k=1}^{\infty}\bigg[\frac{1}{(2n \pi)(k \pi)}
\left\langle \frac{\partial^6 \phi}{\partial x^5 \partial y}, \sqrt{2}\sin(2n  \pi x)\sin(k  \pi y)\right\rangle\\&+\frac{1}{(2n \pi)(k \pi)}\left\langle \frac{\partial^6 \phi}{\partial x \partial y^5}, \sqrt{2}\sin(2n  \pi x)\sin(k  \pi y)\right\rangle\\&+\frac{4C_1C_3}{(2n \pi)(k \pi)} \left\langle\frac{\partial^2 f}{\partial x \partial y}, \sqrt{2}\sin(2n  \pi x)\sin(k  \pi y)\right\rangle\\&+\sum_{n=1}^{\infty}\sum_{k=1}^{\infty}\frac{16C_1(1+\sum_{i=1}^{m}\psi_i\ T^{\alpha-\alpha_i})}{(2n \pi)^2(k \pi)}\bigg\{\bigg|\left\langle \frac{\partial^6 \phi(1,y)}{\partial x^5\partial y}- \frac{\partial^6 \phi(0,y)}{\partial x^5\partial y},\sqrt{2}\sin(k \pi y)\right\rangle\bigg|\\&+\bigg|\left\langle \frac{\partial^6 \phi}{\partial x^5\partial y},\sqrt{2}\cos(2n \pi x)\sin(k \pi y)\right\rangle\bigg|\bigg\}\\&+\sum_{n=1}^{\infty}\sum_{k=1}^{\infty}\frac{8C_1C_3}{(2n \pi)(k \pi)^2}\left\langle f_{xx},\sqrt{2}\sin(2n  \pi x)\cos(k \pi y)\right\rangle.
\end{align*}

The inequality $2ab\leq a^2+b^2$ and the fact $\{\sqrt{2}\sin(2n \pi x) \sin(k \pi y)\}_{n, k=1}^{\infty}$ form orthonormal sequences in $L_2(0,1)$ and $L_2(\Omega)$ respectively, alongside Bessel's inequality, allow us to write
\begin{align*}
 \sum_{n=1}^{\infty}
  \sum_{k=1}^{\infty}\sigma_{nk}|T_{(2n-1)k}(t)|\leq&
2{C_1}(1+\sum_{i=1}^{m}\psi_i\ T^{\alpha-\alpha_i})\biggl\{ \sum_{n=1}^{\infty}\sum_{k=1}^{\infty}\frac{2}{(2n \pi)^2(k \pi)^2}\\&+\bigg\|\frac{\partial^6 \phi}{\partial x^5\partial y}\bigg\|_2^2+\bigg\| \frac{\partial^6 \phi}{\partial x\partial y^5}\bigg\|_2^2\bigg\}+2C_1C_3\sum_{n=1}^{\infty}\sum_{k=1}^{\infty}\bigg\{\frac{1}{(2n \pi)^2(k \pi)^2}+\bigg\| \frac{\partial^2f}{\partial x \partial y}\bigg\|_3^2\bigg\}\\&+8C_1(1+\sum_{i=1}^{m}\psi_i\ T^{\alpha-\alpha_i})\bigg\{\sum_{n=1}^{\infty}\sum_{k=1}^{\infty}\frac{1}{(2n \pi)^4(k \pi)^2}\\&+2\bigg\|\frac{\partial^6 \phi(1,y)}{\partial x^5\partial y}- \frac{\partial^6 \phi(0,y)}{\partial x^5\partial y}\bigg\|_2^2\\&+2\bigg\| \frac{\partial^6 \phi}{\partial x^5\partial y}\bigg\|_2^2\bigg\}+
4C_1C_3\sum_{n=1}^{\infty}\sum_{k=1}^{\infty}\bigg\{\frac{1}{(2n \pi)^2(k \pi)^4}+\bigg\| \frac{\partial^2f}{\partial x^2}\bigg\|_3^2\bigg\}.\end{align*}

\begin{itemize}
\item \textbf{Uniform convergence of series representation of $ \bigl({}^cD_{0_+,t}^{\alpha}+\displaystyle\sum_{i=1}^{m}\psi_i
{}^cD_{0_+,t}^{\alpha_i}\bigr)u(x,y,t)$}
\end{itemize}
By using \eqref{fdesT0}-\eqref{fdesT2nk}, we have following relation
\begin{align*}
 \bigl({}^cD_{0_+,t}^{\alpha}+\sum_{i=1}^{m}\psi_i
{}^cD_{0_+,t}^{\alpha_i}\bigr)u(x,y,t)&=\sum_{k=1}^{\infty}\left(-\mu_kT_{0k}(t)+a(t)f_{0k}(t)\right)Z_{0k}(x,y)\\&+
 \sum_{n=1}^{\infty} \sum_{k=1}^{\infty}\big(-\sigma_{nk}T_{(2n-1)k}(t)+4(\lambda_n)^{3/4}
  T_{2nk}(t)\\&+a(t)f_{(2n-1)k}(t)\big)Z_{(2n-1)k}(x,y)+
\sum_{n=1}^{\infty} \sum_{k=1}^{\infty}\big(-\sigma_{nk}T_{2nk}(t)+a(t)f_{2nk}(t)\big)Z_{2nk}(x,y).
\end{align*}
The series $\displaystyle\sum_{k=1}^{\infty}\mu_kT_{0k}(t),\, \sum_{n=1}^{\infty} \sum_{k=1}^{\infty}\sigma_{nk}T_{(2n-1)k}(t)$ and $ \displaystyle\sum_{n=1}^{\infty} \sum_{k=1}^{\infty}\sigma_{nk}T_{2nk}(t)$ are already proved to be uniformly convergent and using the fact that $|a(t)|\in C[0,T]$ and $|Z_{nk}(x,y)|\leq 4\sqrt{2},\; n,k\in \mathbb{Z}^+$, it remains to show that $\displaystyle\sum_{k=1}^{\infty}f_{0k}(t),\,$ $ \sum_{n=1}^{\infty} \sum_{k=1}^{\infty}f_{(2n-1)k}(t)$ and $ \displaystyle\sum_{n=1}^{\infty} \sum_{k=1}^{\infty}f_{2nk}(t)$ are uniformly convergent.

 Due to Lemma \ref{fest} and Bessel inequality, we can write
\begin{align*}
\bigl|\sum_{k=1}^{\infty}f_{0k}(t)\bigr|&\leq\sum_{k=1}^{\infty}\frac{1}{(k \pi)^2}
\bigl|\left\langle f_{yy},W_{0k}(x,y)\right\rangle\bigr| \leq \sum_{k=1}^{\infty}\frac{4\sqrt{2}}{(k \pi)^2}\big\|f_{yy}\big \|_3,\\
\sum_{n=1}^{\infty} \bigl|\sum_{k=1}^{\infty}f_{(2n-1)k}(t)\bigr|&\leq\sum_{n=1}^{\infty} \sum_{k=1}^{\infty}\frac{4}{(2n \pi)(k \pi)}\bigl|\left\langle f_{xy},\sqrt{2}\sin(2n  \pi x)\sin(k \pi y)\right\rangle\bigr|\\&\leq\sum_{n=1}^{\infty} \sum_{k=1}^{\infty}2\bigg(\frac{1}{(2n \pi)^2(k \pi)^2}+\bigl|\left\langle f_{xy},\sqrt{2}\sin(2n  \pi x)\sin(k \pi y)\right\rangle\bigr|^2\bigg)\\&\leq2\bigg(\sum_{n=1}^{\infty} \sum_{k=1}^{\infty}\frac{1}{(2n \pi)^2(k \pi)^2}+\big\| f_{xy}\big\|^2_2\bigg)\\
\bigl|\sum_{n=1}^{\infty} \sum_{k=1}^{\infty}f_{2nk}(t)\bigr|&\leq\sum_{n=1}^{\infty} \sum_{k=1}^{\infty}\frac{4}{(2n \pi)(k \pi)}\bigl|\left\langle f_{xy},\sqrt{2}\cos(2n  \pi x)\sin(k \pi y)\right\rangle\bigr| \\&\leq\sum_{n=1}^{\infty} \sum_{k=1}^{\infty}2\bigg(\frac{1}{(2n \pi)^2(k \pi)^2}+\bigl|\left\langle f_{xy},\sqrt{2}\cos(2n  \pi x)\sin(k \pi y)\right\rangle\bigr|^2\bigg)\\&\leq 2\left(\sum_{n=1}^{\infty} \sum_{k=1}^{\infty}\frac{1}{(2n \pi)^2(k \pi)^2}+\big\| f_{xy}\big\|^2_2\right).
\end{align*}
Hence, $ \bigl({}^cD_{0_+,t}^{\alpha}+\displaystyle\sum_{i=1}^{m}\psi_i
{}^cD_{0_+,t}^{\alpha_i}\bigr)u(x,t)$ represents a continuous function.

After showing the regularity of the obtained solution, we will show the uniqueness of the solution.

\textbf{Uniqueness of the solution:}\\
In order to establish the uniqueness of the inverse source problem \eqref{problem}-\eqref{odc}, we need to prove the uniqueness of time dependent source term $a(t)$ and $u(x,y,t).$

The expression of  $a(t)$, given by  \eqref{a(t)}, involves definite integral of given function $f(x,y,t)$ and fractional derivative of known function $E(t).$ One-one nature of operators (definite integral and fractional derivative) ensures the uniqueness of $a(t).$

For uniqueness of $u(x,y,t),$ let us consider $u(x,y,t)$ and $v(x,y,t)$ be two solutions, and let
$\bar{u}(x,y,t)=u(x,y,t)-v(x,y,t)$.
Then $\bar{u}(x,y,t)$ satisfy the  equation
\begin{equation}\label{problemunique}
\bigl({}^cD_{0_+,t}^{\alpha}+\sum_{i=1}^{m}\psi_i
{}^cD_{0_+,t}^{\alpha_i}\bigr)\bar{u}(x,y,t)+\Delta^2\bar{u}(x,y,t)=0,\qquad (x,y,t)\in \Pi:=\Omega \times (0,T).
\end{equation}
Dirichlet boundary condition
\begin{eqnarray}\label{dbcsu}
&\bar{u}_{x}(0,y,t)=0=\bar{u}_{xxx}(0,y,t),\;
\bar{u}_{y}(x,0,t)=0=\bar{u}_{y}(x,1,t),&\;\nonumber\\
&\bar{u}_{yyy}(x,0,t)=0=\bar{u}_{yyy}(x,1,t),&
\end{eqnarray}
nonlocal boundary conditions
\begin{equation}\label{nlbcsu}
\bar{u}(0,y,t)=\bar{u}(1,y,t),\; \bar{u}_{xx}(0,y,t)=\bar{u}_{xx}(1,y,t),
\end{equation}
and initial condition
\begin{equation}\label{icu}
\bar{u}(x,y,0)=0,\qquad (x,y)\in \Omega:(0,1)\times(0,1),
\end{equation}
Consider the functions
\begin{gather*} 
\bar {T}_{0k}(t)=\int_0^{1}\int_0^{1}\bar u(x,y,t)W_{0k}(x,y)dxdy, \\
\bar {T}_{(2n-1)k}(t)=\int_0^{1}\int_0^{1}\bar {u}(x,y,t)W_{(2n-1)k}(x,y)dxdy,\\
\bar {T}_{2nk}(t)=\int_0^{1}\int_0^{1}\bar {u}(x,y,t)W_{2nk}(x,y)dxdy.
\end{gather*}
Taking the fractional derivative, we get
\begin{gather*}
 \bigl({}^cD_{0_+,t}^{\alpha}+\sum_{i=1}^{m}\psi_i
{}^cD_{0_+,t}^{\alpha_i}\bigr)\bar {T}_{0k}(t)=-\mu_k\bar {T}_{0k}(t),  \\
 \bigl({}^cD_{0_+,t}^{\alpha}+\sum_{i=1}^{m}\psi_i
{}^cD_{0_+,t}^{\alpha_i}\bigr)\bar {T}_{(2n-1)k}(t)
  =-\sigma_{nk}\bar {T}_{(2n-1)k}(t)+4(\lambda_n)^{3/4}
  \bar {T}_{2nk}(t),\\
 \bigl({}^cD_{0_+,t}^{\alpha}+\sum_{i=1}^{m}\psi_i
{}^cD_{0_+,t}^{\alpha_i}\bigr)\bar {T}_{2nk}(t) =-\sigma_{nk}\bar {T}_{2nk}(t).
  \end{gather*}
The solution of the above system is
\begin{align*}
\bar{T}_{0k}(t)=&\mbox{\Large $e$ \hspace{-.3cm} }_{(\mu_1(\alpha-\alpha_1),...,\mu_m(\alpha-\alpha_m),
\mu_k\alpha),1}(t)\;\bar{T}_{0k}(0)+
\sum_{i=1}^{m}\psi_i\mbox{\Large $e$ \hspace{-.3cm} }_{(\mu_1(\alpha-\alpha_1),...,
\mu_m(\alpha-\alpha_m),\mu_k\alpha_i),\alpha+1-\alpha_i}(t)\;\bar{T}_{0k}(0)\\
\bar{T}_{2nk}(t)=&\mbox{\Large $e$ \hspace{-.3cm} }_{(\mu_1(\alpha-\alpha_1),...,\mu_m(\alpha-\alpha_m),\sigma_{nk}
\alpha),1}(t)\;\bar{T}_{2nk}(0)
+\sum_{i=1}^{m}\psi_i\mbox{\Large $e$ \hspace{-.3cm} }_{(\mu_1(\alpha-\alpha_1),...,
\mu_m(\alpha-\alpha_m),\sigma_{nk}\alpha_i),\alpha+1-\alpha_i}(t)\;\bar{T}_
{2nk}(0)\\
\bar{T}_{(2n-1)k}(t)=&\mbox{\Large $e$ \hspace{-.3cm} }_{(\mu_1(\alpha-\alpha_1),...,\mu_m(\alpha-\alpha_m),\sigma_{nk}
\alpha),1}(t)\;\bar{T}_{(2n-1)k}(0)+4(\lambda_n)^{3/4}\bar{T}_{2nk}(t)*
\mbox{\Large $e$ \hspace{-.3cm} }_{(\mu_1(\alpha-\alpha_1),...,\mu_m(\alpha-\alpha_m),
\sigma_{nk}\alpha_i),\alpha}(t)\\&+\sum_{i=1}^{m}\psi_i
\mbox{\Large $e$ \hspace{-.3cm} }_{(\mu_1(\alpha-\alpha_1),...,
\mu_m(\alpha-\alpha_m),\sigma_{nk}\alpha_i),\alpha+1-\alpha_i}(t)\;
\bar{T}_{(2n-1)k}(0).
\end{align*}
By using the initial condition \eqref{icu}, we have
$$
\bar {T}_{0k}(t) = 0,\quad \bar {T}_{(2n-1)k}(t)=0,\quad \bar {T}_{2nk}(t)=0,\quad
t\in(0,T).
$$
Consequently, the uniqueness of the solution follows from the completeness
 of the set of function $$\{W_{0k}(x,y), W_{(2n-1)k}(x,y),W_{2nk}(x,y)\},\quad n,\; k\in\mathbb{Z}^+.$$
 \textbf{Stability of the solution:}\\
Now, we will show that the solution depends continuously on given datum, assume that $\{u(x,y,t),a(t)\}$ and $\{\tilde{u}(x,y,t),\tilde{a}(t)\}$ be two solution sets of the inverse source problem (\ref{problem})-(\ref{odc}), corresponding to given data $\{f(x,y,t),\phi(x), E(t)\}$ and $\{\tilde{f}(x,y,t),\tilde \phi(x),\tilde E(t)\},$ respectively.

From \eqref{a(t)}, we have
\begin{eqnarray*}
a(t)-\tilde{a}(t)=&\left(\displaystyle\int_{0}^{1}\int_{0}^{1}f(x,y,t)
dxdy\right)^{-1}
\bigl({}^cD_{0_+,t}^{\alpha}+\displaystyle\sum_{i=1}^{m}\psi_i
{}^cD_{0_+,t}^{\alpha_i}\bigr)E(t)\\&-\left(\displaystyle
\int_{0}^{1}\int_{0}^{1}\tilde{f}(x,y,t)
dxdy\right)^{-1}
\bigl({}^cD_{0_+,t}^{\alpha}+\displaystyle\sum_{i=1}^{m}\psi_i
{}^cD_{0_+,t}^{\alpha_i}\bigr)\tilde{E}(t),
\end{eqnarray*}
which leads us to the relation
\begin{eqnarray}
a(t)-\tilde{a}(t)=&
\frac{\left(
\displaystyle\int_{0}^{1}\int_{0}^{1}\tilde{f}(x,y,t)-f(x,y,t)
dxdy\right)\bigl(\displaystyle{}^cD_{0_+,t}^{\alpha}+\displaystyle\sum_{i=1}^{m}\psi_i
{}^cD_{0_+,t}^{\alpha_i}\bigr)E(t)}{\displaystyle\int_{0}^{1}\int_{0}^{1}f(x,y,t)dxdy
\displaystyle\int_{0}^{1}\int_{0}^{1}\tilde{f}(x,y,t)dxdy}\nonumber\\&+
\frac{\left(\displaystyle\int_{0}^{1}\int_{0}^{1}f(x,y,t)dxdy\right)\bigl(\displaystyle{}^cD_{0_+,t}^{\alpha}+\displaystyle\sum_{i=1}^{m}\psi_i
{}^cD_{0_+,t}^{\alpha_i}\bigr)
(E(t)-\tilde{E(t)})}{\displaystyle\int_{0}^{1}\int_{0}^{1}f(x,y,t)dxdy
\displaystyle\int_{0}^{1}\int_{0}^{1}\tilde{f}(x,y,t)dxdy}\label{a(t)stability1}.
\end{eqnarray}
Note that there exist constants $C_4,\, C_5,\, C_6$ and $C_7$ such that
\begin{eqnarray*}
\|\bigl({}^cD_{0_+,t}^{\alpha}+\sum_{i=1}^{m}\psi_i
{}^cD_{0_+,t}^{\alpha_i}\bigr)E(t)\|_1\leq C_4,\\
\|\bigl({}^cD_{0_+,t}^{\alpha}+\sum_{i=1}^{m}\psi_i
{}^cD_{0_+,t}^{\alpha_i}\bigr)(E(t)-\tilde{E}(t)\|_1\leq C_5\|E-\tilde{E}\|_1,\\
\bigg\|\displaystyle\int_{0}^{1}\int_{0}^{1}f(x,y,t)
dxdy\bigg\|_1\leq C_6,\\
\bigg\|\displaystyle\int_{0}^{1}\int_{0}^{1}\tilde{f}(x,y,t)-f(x,y,t)
dxdy\bigg\|_1\leq C_7\|f-\tilde{f}\|_3.\end{eqnarray*}

Consequently, \eqref{a(t)stability1} becomes
\begin{eqnarray*}\label{a(t)stability2}
\|a-\tilde{a}\|_1\leq C_2^2(C_4C_7\big\|f-\tilde{f}\big\|_3+C_5C_6\big\|E-\tilde{E}\big\|_1).
\end{eqnarray*}
After establishing the stability result for source term $a(t)$, we will establish the inequality of $u(x,y,t)$ that ensures the dependence of $u(x,y,t)$ on the given datum.

From \eqref{seriesu}, we have
\begin{align*}
|u(x,y,t)-\tilde{u}(x,y,t)|&\leq\displaystyle\sum_{k=1}^{\infty}|T_{0k}(t)-\tilde{T}_{0k}(t)|
Z_{0k}(x,y)\nonumber+
  \displaystyle\sum_{n=1}^{\infty}
  \displaystyle\sum_{k=1}^{\infty}\big[|T_{(2n-1)k}(t)-\tilde{T}_{(2n-1)k}(t)|
  Z_{(2n-1)k}(x,y)\nonumber\\&
  +|T_{2nk}(t)-\tilde{T}_{2nk}(t)Z_{2nk}(x,y)\big],
\end{align*}
since $|Z_{nk}(x,y)|\leq\; \sqrt{2},\quad \forall n\in \mathbb{Z}^+_0,\;k\in\mathbb{Z}^+$, so we have
\begin{align}
|u(x,y,t)-\tilde{u}(x,y,t)|&\leq\sqrt{2}\sum_{k=1}^{\infty}|T_{0k}(t)-\tilde{T}_{0k}(t)|+
  \sqrt{2}\sum_{n=1}^{\infty}
  \sum_{k=1}^{\infty}\big[|T_{(2n-1)k}(t)-\tilde{T}_{(2n-1)k}(t)| \nonumber\\&+|T_{2nk}-\tilde{T}_{2nk}(t)\big]\label{ustability},
\end{align}
By using \eqref{T0kest}, we have
\begin{align}
\bigg|\sum_{k=1}^{\infty}(T_{0k}(t)-\tilde{T}_{0k}(t))\bigg|\leq & \sum_{k=0}^{\infty}\frac{2\sqrt{2}C_1\big\|\phi-\tilde{\phi}\big\|_2}{(k \pi)^4}\biggl(\frac{1}{t^\alpha}+
\sum_{i=1}^{m}\frac{\psi_i}{t^{\alpha_i}}\biggr)+
\sum_{k=0}^{\infty}\frac{2\sqrt{2}C_1C_3\big\|f-\tilde{f}\big\|_3}{(k \pi)^4}
\label{T0kstability}.
\end{align}
Inequality \eqref{T2nkest} is used to get the following
\begin{align}
\bigg|\sum_{n=1}^{\infty}\sum_{k=1}^{\infty}(T_{2nk}(t)-\tilde{T}_{2nk}(t))
\bigg|\leq & \sum_{n=1}^{\infty}\sum_{k=1}^{\infty}\frac{C_1\big\|\phi-\tilde{\phi}\big\|_2}
{\sqrt{2}(n \pi)^2(k \pi)^2}\biggl(\frac{1}{t^\alpha}+
\sum_{i=1}^{m}\frac{\psi_i}{t^{\alpha_i}}\biggr)
\nonumber\\&+\sum_{n=1}^{\infty}\sum_{k=1}^{\infty}\frac{C_1C_3\big\|f-\tilde{f}\big\|_3}{\sqrt{2}(n \pi)^2
(k \pi)^2}\label{T2nkstability}.
\end{align}
Lastly for $\displaystyle\sum_{n=1}^{\infty}\sum_{k=1}^{\infty}(T_{(2n-1)k}-
\tilde{T}_{(2n-1)k}(t))$,
we will make use of \eqref{T2nkest}
\begin{align}
\bigg|\sum_{n=1}^{\infty}\sum_{k=1}^{\infty}
(T_{(2n-1)k}(t)-\tilde{T}_{(2n-1)k}(t))\bigg|\leq &\sum_{n=1}^{\infty}\sum_{k=1}^{\infty}\frac{C_1\|\phi-\tilde{\phi}\|_2}{\sqrt{2}(n \pi)^2(k \pi)^2}\biggl(\frac{1}{t^\alpha}+
\sum_{i=1}^{m}\frac{\psi_i}{t^{\alpha_i}}\biggr)+
\sum_{n=1}^{\infty}\sum_{k=1}^{\infty}\frac{C_1C_3\|f-\tilde{f}\|_3}{\sqrt{2}
(n \pi)^2(k \pi)^2}\nonumber
\\&+\sum_{n=1}^{\infty}\sum_{k=1}^{\infty}\frac{C^2_1\|\phi-\tilde{\phi}\|_2}{2(n \pi)^3(k \pi)^4}\biggl(\frac{1}{t^\alpha}+
\sum_{i=1}^{m}\frac{\psi_i}{t^{\alpha_i}}\biggr)\nonumber\\&
+\sum_{n=1}^{\infty}\sum_{k=1}^{\infty}\frac{C^2_1C_3\| f-\tilde{f}\|_3}{2(n \pi)^3(k \pi)^4}\label{T2n1kstability}.
\end{align}
By using \eqref{T0kstability}-\eqref{T2n1kstability} in \eqref{ustability} takes the form
\begin{align*}
|u(x,y,t)-\tilde{u}(x,y,t)|\leq&\displaystyle\sum_{k=0}^{\infty}
\frac{4C_1\|\phi-\tilde{\phi}\|_2}{(k \pi)^4}\biggl(\frac{1}{t^\alpha}+
\displaystyle\sum_{i=1}^{m}\frac{\psi_i}{t^{\alpha_i}}\biggr)+
\sum_{k=0}^{\infty}\frac{4C_1C_3\|f-\tilde{f}\|_3}{(k \pi)^4}\nonumber\\
  &+\displaystyle\sum_{n=1}^{\infty}\sum_{k=1}^{\infty}
  \frac{2C_1\|\phi-\tilde{\phi}\|_2}
{(n \pi)^2(k \pi)^2}\biggl(\frac{1}{t^\alpha}+
\displaystyle\sum_{i=1}^{m}\frac{\psi_i}{t^{\alpha_i}}\biggr)
+\displaystyle\sum_{n=1}^{\infty}\sum_{k=1}^{\infty}\frac{2C_1C_3\|f-\tilde{f}\|_3}{(n \pi)^2(k \pi)^2}\nonumber
\\&+\sum_{n=1}^{\infty}\sum_{k=1}^{\infty}\frac{C^2_1\|\phi-\tilde{\phi}\|_2}{\sqrt{2}(n \pi)^3(k \pi)^4}\biggl(\frac{1}{t^\alpha}+
\sum_{i=1}^{m}\frac{\psi_i}{t^{\alpha_i}}\biggr)
+\sum_{n=1}^{\infty}\sum_{k=1}^{\infty}\frac{C^2_1C_3\| f-\tilde{f}\|_3}{\sqrt{2}(n \pi)^3(k \pi)^4}.
\label{ustability1}
\end{align*}
Since, the double series $\displaystyle\sum_{n=1}^{\infty}\sum_{k=1}^{\infty}\frac{1}{n^2k^2}$ and $\displaystyle\sum_{n=1}^{\infty}\sum_{k=1}^{\infty}\frac{1}{n^3k^4}$ alongside the series $\displaystyle\sum_{k=1}\frac{1}{k^4}$  are convergent, so there exists constants $C_8$ and $C_9$ such that
\begin{align*}
\|u-\tilde{u}\|_3\leq C_8\|\phi-\tilde{\phi}\|_2+C_9\|f-\tilde{f}\|_3,
\end{align*}
where $C_8$ and $C_9$ are constants independent of $n$ and $k$.
\end{proof}

\end{document}